\documentclass[11pt]{amsart}
\usepackage{amsmath}
\usepackage{amssymb}
\usepackage{latexsym}
\usepackage{graphicx}

\newcommand{\N}{\mathbb{N}}

\newtheorem{theorem}{Theorem}
\newtheorem{corollary}{Corollary}
\newtheorem{lemma}{Lemma}
\theoremstyle{remark}
\newtheorem{remark}{Remark}

\begin{document}

\title[integral means by iteration of Janowski functions]{\large integral means of analytic mappings by iteration of Janowski functions}

\author[K. O. Babalola]{K. O. BABALOLA}

\begin{abstract}
In this short note we apply certain iteration of the Janowski functions to estimate the integral means of some analytic and univalent mappings of $|z|<1$. Our method of proof follows an earlier one due to Leung \cite{YJ}.
\end{abstract}



\maketitle

\section{Introduction}

Let $A$ be the class of normalized analytic functions $f(z)=z+a_2z^2+...$ in the unit disk $|z|<1$. In \cite{KO}, among others, we added a new generalization class, namely; $T_n^\alpha[a,b]$, $\alpha>0$, $-1\leq b<a\leq 1$ and $n\in \N$; to the large body of analytic and univalent mappings of the unit disk $|z|<1$. This consists of functions in $|z|<1$ satisfying the geometric conditions
$$\frac{D^nf(z)^\alpha}{\alpha^nz^\alpha}\in P[a,b]\eqno{(1)}$$
where $P[a,b]$ is the family of Janowski functions $p(z)=1+c_1z+\cdots$ which are subordinate to $L_0(a,b:z)=(1+az)/(1+bz)$, $-1\leq b<a\leq 1$, in $|z|<1$. The operator $D^n$, defined as $D^nf(z)=z[D^{n-1)}f(z)]'$ with $D^0f(z)=f(z)$, is the well known Salagean derivative \cite{GS}.

In Section 2 of the paper \cite{KO} we extended certain integral iteration of the class of Caratheodory functions (which we developed in \cite{BO}) to $P[a,b]$ via which the new class, $T_n^\alpha[a,b]$, was studied. The extension was obtained simply by choosing the analytic function $p(z)=1+c_1z+\cdots$, Re $p(z)>0$ from $P[a,b]$ in the iteration defined in \cite{BO} as:
$$p_n(z)=\frac{\alpha}{z^\alpha}\int_0^zt^{\alpha-1}p_{n-1}(t)dt,\;\;n\geq 1,$$
with $p_0(z)=p(z)$.

We will denote this extention by $P_n[a,b]$ in this note. We had remarked (in \cite{KO}) that the statements (i) $p(z)\prec L_0(a,b:z)$, (ii) $p\in P[a,b]$, (iii) $p_n(z)\in P_n[a,b]$ and (iv) $p_n\prec L_n(a,b:z)$ are all equivalent. Thus we also remarked that (1) is equivalent to $f(z)^\alpha/z^\alpha\in P_n[a,b]$. This new equivalent geometric condition will lead us to the following interesting results regarding the integral means of functions in $T_n^\alpha[a,b]$ for $0<\alpha\leq 1$ and $n\geq 1$.

\begin{theorem}
Let $\Phi$ be a convex non-decreasing function $\Phi$ on $(-\infty$, $\infty)$. Then for $f\in T_n^\alpha[a,b]$, $\alpha\in(0,1]$, $n\geq 1$ and $r\in(0,1)$
$$\int_{-\pi}^\pi\Phi\left(\log|f^\prime(re^{i\theta})|\right)d\theta\leq\int_{-\pi}^\pi\Phi\left(\log\left|\frac{L_{n-1}(a,b:re^{i\theta})k'(re^{i\theta})^{1-\alpha}}{L_0(re^{i\theta})^{1-\alpha}}\right|\right)d\theta\eqno{(2)}$$
where
$$L_n(a,b:z)=\frac{\alpha}{z^\alpha}\int_0^zt^{\alpha-1}L_{n-1}(a,b:t)dt,\;\;\;n\geq 1$$
and $k(z)=z/(1-z)^2$ is the Koebe function.
\end{theorem}

\begin{theorem}
With the same hypothesis as in Theorem $1$, we have 
$$\int_{-\pi}^\pi\Phi\left(-\log|f^\prime(re^{i\theta})|\right)d\theta\leq\int_{-\pi}^\pi\Phi\left(-\log\left|\frac{L_{n-1}(a,b:re^{i\theta})k'(re^{i\theta})^{1-\alpha}}{L_0(re^{i\theta})^{1-\alpha}}\right|\right)d\theta.$$
\end{theorem}

The above inequalities represent the integral means of functions of the class $T_n^\alpha[a,b]$ for $\alpha\in(0,1]$ and $n\geq 1$. Our method of proof follows an earlier one due to Leung \cite{YJ} using the equivalent geometric relations $f(z)^\alpha/z^\alpha\in P_n[a,b]$ for $f\in T_n^\alpha[a,b]$.\vskip 2mm

It is worthy of note that very many particular cases of the above results can be obtained by specifying the parameters $n,\alpha,a$ and $b$ as appropriate. In particular, the following  special cases of $P[a,b]$ are well known: $P[1,-1]$; $P[1-2\beta,-1]$, $0\leq\beta<1$; $P[1,1/\beta-1]$, $\beta>1/2$; $P[\beta,-\beta]$, $0<\beta\leq 1$ and $P[\beta,0]$, $0<\beta\leq 1$ (see \cite{KO}). Thus several cases of $T_n^\alpha[a,b]$ may also be deduced.

\section{Fundamental Lemmas}

\medskip

The following results are due to Baernstein \cite{AB} and Leung \cite{YJ}. Let $g(x)$ be a real-valued integrable function on $[-\pi,\pi]$. Define $g^\ast(x)=\sup_{|E|=2\theta}\int_Eg$, ($0\leq\theta\leq\pi$) where $|E|$ denotes the Lebesgue measure of the set $E$ in $[-\pi,\pi]$. Further details can be found in the Baernstein's work \cite{AB}.

\begin{lemma}[\cite{AB}]
For $g,h\in L^1[-\pi,\pi]$, the following statements are equivalent:

{\rm (i)} For every convex non-decreasing function $\Phi$ on $(-\infty$, $\infty)$,
$$\int_{-\pi}^\pi\Phi(g(x))dx\leq\int_{-\pi}^\pi\Phi(h(x))dx.$$

{\rm (ii)} For every $t\in(-\infty$, $\infty)$,
$$\int_{-\pi}^\pi[g(x)-t]^+dx\leq\int_{-\pi}^\pi[h(x)-t]^+dx.$$

{\rm (iii)} $g^\ast(\theta)\leq h^\ast(\theta)$, $(0\leq\theta\leq\pi)$.
\end{lemma}

\begin{lemma}[\cite{AB}]
If $f$ is normalized and univalent in $|z|<1$, then for each $r\in(0,1)$, $(\pm\log|f(re^{i\theta}|)^\ast\leq(\pm\log|k(re^{i\theta}|)^\ast$.
\end{lemma}

\begin{lemma}[\cite{YJ}]
For $g,h\in L^1[-\pi,\pi]$, $[g(\theta)+h(\theta)]^\ast\leq g^\ast(\theta)+h^\ast(\theta)$. Equality holds if $g,h$ are both symmetric in $[-\pi,\pi]$ and nonincreasing in $[0,\pi]$.
\end{lemma}

\begin{lemma}[\cite{YJ}]
If $g,h$ are subharmonic in $|z|<1$ and $g$ is subordinate to $h$, then for each $r\in(0,1)$,  $g^\ast(re^{i\theta})\leq h^\ast(re^{i\theta})$, $(0\leq\theta\leq\pi)$.
\end{lemma}

\begin{corollary}
If $p\in P_n[a,b]$, then
$$\left(\pm\log|p_n(re^{i\theta}|\right)^\ast\leq\left(\pm\log|L_n(a,b:re^{i\theta}|\right)^\ast,\;\;0\leq\theta\leq\pi.$$
\end{corollary}

\begin{proof}
Since $p_n(z)$ and $L_n(a,b:z)$ are analytic, $\log|p_n(z)|$ and $\log|L_n(a,b:z)|$ are both subharmonic in $|z|<1$. Furthermore, since $p_n\prec L_n(a,b:z)$, there exists $w(z)$ ($|w(z)|<1$), such that $p_n(z)=L_n(a,b:w(z))$. Thus we have $\log p_n(z)=\log L_n(a,b:w(z))$ so that $\log p_n(z)\prec\log L_n(a,b:z)$. Hence by Lemma 3 we have the first of the inequalities.

As for the second, we also note from the above that $1/p_n(z)=1/L_n(a,b:w(z))$ so that $-\log p_n(z)=-\log L_n(a,b:w(z))$ and thus $-\log p_n(z)\prec-\log L_n(a,b:z)$. Also $\log|1/p_n(z)|$ and $\log|1/L_n(a,b:z)|$ are both subharmonic in $|z|<1$ since $1/p_n(z)$ and $1/L_n(a,b:z)$ are analytic there. Thus by Lemma 3 again, we have the desired inequality.
\end{proof}

\section{Proofs of Main Results}
 \medskip

We begin with

\begin{proof}[Proof of Theorem $1$]
Since $f\in T_n^\alpha[a,b], \alpha\in(0,1]$, then there exists $p_n\in P_n[a,b]$, such that $f(z)^\alpha/z^\alpha=p_n(z)$. Then $f'(z)=p_{n-1}(z)(f(z)/z)^{1-\alpha}$ so that
$$\aligned
\log|f'(z)|
&=\log|p_{n-1}(z)|+\log\left|\frac{f(z)}{z}\right|^{1-\alpha}\\
&=\log|p_{n-1}(z)|+(1-\alpha)\log\left|\frac{f(z)}{z}\right|
\endaligned\eqno{(3)}$$
so that, by Lemma 3,
$$(\log|f'(z)|)^\ast=(\log|p_{n-1}(z)|)^\ast+\left(\log\left|\frac{f(z)}{z}\right|^{1-\alpha}\right)^\ast.$$
For $n\geq 1$, $f(z)$ is univalent (see \cite{KO}), so that by Lemma 2 and Corollary 1 we have
$$\aligned
(\log|f'(z)|)^\ast
&=\left(\log|L_{n-1}(a,b:re^{i\theta}|\right)^\ast+\left(\log\left|\frac{k(re^{i\theta})}{r}\right|^{1-\alpha}\right)^\ast\\
&=\left(\log\left|\frac{L_{n-1}(a,b:re^{i\theta})k'(re^{i\theta})^{1-\alpha}}{L_0(re^{i\theta})^{1-\alpha}}\right|\right)^\ast.
\endaligned$$
Hence by Lemma 1, we have the inequality. If for some $r\in(0,1)$ and some strictly convex $\Phi$, we consider the function $f_0(z)$ is defined by
$$e^{-i\alpha\gamma}\frac{f_0(ze^{i\gamma})^\alpha}{z^\alpha}=L_n(a,b:ze^{i\gamma})\eqno{(4)}$$
for some real $\gamma$. Then we have
$$\aligned
e^{i\gamma(1-\alpha)}\frac{f_0(ze^{i\gamma})^{\alpha-1}f_0'(ze^{i\gamma})}{z^{\alpha-1}}&=L_n(a,b:ze^{i\gamma})+\frac{ze^{i\gamma}L_n(a,b:ze^{i\gamma})}{\alpha}\\
&=L_{n-1}(a,b:ze^{i\gamma}),\endaligned$$
so that
$$|f_0'(ze^{i\gamma})|=|L_{n-1}(a,b:ze^{i\gamma})|\left|\frac{f_0(ze^{i\gamma})}{ze^{i\gamma}}\right|^{1-\alpha}.$$
Now equality in (2) can be attained by taking $|f_0(z)|=|k(z)|$. This completes the proof.
\end{proof}

Next we have

\begin{proof}[Proof of Theorem $2$]
From (3) we have
$$\log\frac{1}{|f'(z)|}=\log\frac{1}{|p_{n-1}(z)|}+(1-\alpha)\log\left|\frac{z}{f(z)}\right|.$$
Hence, by Lemmas 2, 3 and Corollary 1 again, we have
$$\aligned
(-\log|f'(z)|)^\ast
&\leq\left(\log\left|\frac{1}{L_{n-1}(a,b:re^{i\theta})}\right|\right)^\ast+\left(\log\left|\frac{r}{k(re^{i\theta})}\right|^{1-\alpha}\right)^\ast\\
&=\left(-\log\left|\frac{L_{n-1}(a,b:re^{i\theta})k'(re^{i\theta})^{1-\alpha}}{L_0(re^{i\theta})^{1-\alpha}}\right|\right)^\ast.
\endaligned$$
Hence by Lemma 1, we have the inequality. Similarly if equality is attained for some $r\in(0,1)$ and some strictly convex $\Phi$, then $f_0(z)$ given by (4) is the equality function.
\end{proof}

\section{Particular cases}
 \medskip

With the same hypothesis as in Theorem $1$ except:
 \medskip
 
(i) $n=1$, we have:
 
$$\int_{-\pi}^\pi\Phi\left(\pm\log|f'(re^{i\theta})|\right)d\theta\leq\int_{-\pi}^\pi\Phi\left(\pm\log\left|\frac{L_0(a,b:re^{i\theta})k'(re^{i\theta})^{1-\alpha}}{L_0(re^{i\theta})^{1-\alpha}}\right|\right)d\theta.$$
 \medskip

(ii) $\alpha=1$, we have:
 
$$\int_{-\pi}^\pi\Phi\left(\pm\log|f'(re^{i\theta})|\right)d\theta\leq\int_{-\pi}^\pi\Phi\left(\pm\log|L_{n-1}(a,b:re^{i\theta})|\right)d\theta.$$
 \medskip

(iii) $n=\alpha=1$, we have:
 
$$\int_{-\pi}^\pi\Phi\left(\pm\log|f^\prime(re^{i\theta})|\right)d\theta\leq\int_{-\pi}^\pi\Phi\left(\pm\log|L_0(a,b:re^{i\theta})|\right)d\theta.$$
 \medskip

\begin{remark}
The case $n=1, a=1$ and $b=-1$ gives the estimate for the special case $s(z)=z$ of the Leung results \cite{YJ}.
\end{remark}
 \medskip

{\it Acknowledgements.} The author acknowledges the Abdus Salam International Centre for Theoretical Physics, Trieste, Italy for providing the research paper \cite{YJ}.

\vspace{10pt}

\hspace{-4mm}{\small{Received}}

\vspace{-12pt}
\ \hfill \
\begin{tabular}{c}
{\small\em  Department of Mathematics}\\
{\small\em  University of Ilorin}\\
{\small\em  Ilorin, Nigeria}\\
{\small\em E-mail: {\tt kobabalola@gmail.com}} \\
\end{tabular}

\end{document}